\documentclass{article}

 \pdfoutput=1

\usepackage{graphicx}
\usepackage{amsmath}
\usepackage{amsfonts}
\usepackage{amssymb}
\usepackage{amsthm}
\usepackage[all]{xy}
\usepackage{fancyhdr}
\usepackage{subfig}
\usepackage{url}
\usepackage{cite} 
\usepackage{underscore}
\usepackage{appendix}

\usepackage{epsf}

\theoremstyle{definition}
\newtheorem{theorem}{Theorem}[section]

\newtheorem{defn}{Definition}[section]
\newtheorem{lemma}{Lemma}[section]

\newtheorem{note}{Note}[section]

\newcommand{\R}{\mathbb{R}}

\newcommand{\N}{\mathbb{N}}

\newcommand{\ve}{\mathbf} 

\begin{document}

\title{A Max-Plus Model of Asynchronous Cellular Automata}
\author{
Ebrahim~L.~Patel\footnote{Corresponding author: ebrahim.patel@maths.ox.ac.uk} \\
{\small Mathematical Institute, University of Oxford } \\ 
{\small Oxford OX2 6GG, UK } \\  
and \\
David Broomhead  \\
{\small School of Mathematics, University of Manchester } \\
{\small Manchester M13 9PL, UK }
} 

\maketitle    

\begin{abstract}
This paper presents a new framework for asynchrony.   This has its origins in our attempts to better harness the internal decision making process of cellular automata (CA).  Thus, we show that a max-plus algebraic model of asynchrony arises naturally from the CA requirement that a
cell receives the state of each neighbour before updating.  The
significant result is the existence of a bijective mapping between
the asynchronous system and the
synchronous system classically used to update cellular automata.  Consequently, although the CA outputs look qualitatively different, when surveyed on ``contours" of real time, the asynchronous CA replicates the synchronous CA.  Moreover, this type of asynchrony
is simple - it is characterised by the underlying network structure of the cells, and long-term behaviour is deterministic and periodic due to the
linearity of max-plus algebra.  The findings lead us to proffer max-plus algebra as: (i) a more accurate and efficient underlying timing mechanism for models of patterns seen in nature, and (ii)  a foundation for promising extensions and  applications.
%
\end{abstract}

\section{Introduction}
\begin{figure}[!hbp]
\centering
  \includegraphics[scale=0.4,trim=0 0 0 0,clip]{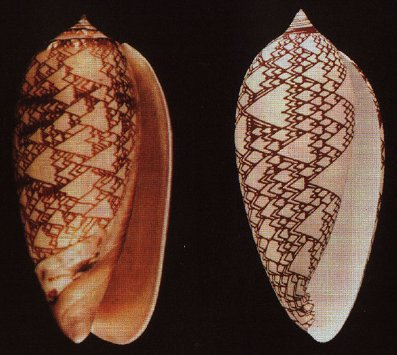}
  \caption{Seashell patterns: left is real, right generated by a CA.  Source: \cite{meinhardt} }
  \label{fig:seashell}
\end{figure}

A cellular automaton (or CA, where we also abbreviate the plural ``cellular automata" to CA, allowing the context to save confusion) is a discrete dynamical system, consisting of an array of identical cells, each possessing a state.  The states evolve, according to some local rule, in discrete time steps.  The first CA models were synchronous, where all cells update once on each time step, and the foundations of the study of these CA were laid by Wolfram in the 1980s \cite{{Wol1},{Wol2}}.  A popular application of such CA is the construction of models of natural growth processes such as seashell patterns and snowflakes \cite{meinhardt, snowflake}.  Figure~\ref{fig:seashell} (right) shows a CA pattern typically examined by Wolfram inscribed on a seashell; the similarity to the real seashell pattern (on the left) is evident.
A natural extension is the introduction of asynchronous update times. Indeed, in terms of seashell patterns, Gunji demonstrated that different asynchronous update methods yield different patterns, leading to their
conjecture that asynchrony is intrinsic to living systems \cite{Gunji}.

A preliminary observation of asynchronous cellular automata was
made in \cite{Ing}, where the authors compared the properties of
synchronous CA with two types of CA that iterate asynchronously.  Subsequent studies employed the methods of \cite{Ing} as special
cases to conduct specific studies into asynchronous CA \cite{Bersini, Fates1, Fates2, LeCaer, Sch}.  Many of these authors attest to asynchrony as being stochastic in nature. This is a general viewpoint in
light of their applications: such asynchrony relies on continuous time \cite{Sch} and is also
likely to be more robust \cite{Sch, Bersini, Fates2}, thereby aiding a better description of biological
phenomena. For example, given a system of coupled cells, the update of cell states
depends on a predefined probability \cite{Ing, Sch}. This consequently also led
Sch\"{o}nfisch and de Roos to conjecture that, while synchronous updating can produce
periodic orbits, asynchronous systems will only yield patterns
that converge to a fixed point or patterns that are chaotic \cite{Sch}.

The argument for asynchronous updating being stochastic has been challenged by authors such as Cornforth et al in \cite{Cor2}. They claim that such probabilistic updating schemes are used because of the
oversimplification of biologically inspired models. They further argue for mimicking
appropriate aspects of nature more closely to create better computational models.
Thus, the authors have drawn
attention to a large class of behaviours of natural processes, in which the
updating is asynchronous but not stochastic \cite{Cor1}. 

Moreover, underlying a synchronous update scheme is the notion of a `global clock', in the sense that it assigns the same update time to all cells.  In recent years, the disadvantages of synchrony in this context have been exposed \cite{fates2013}.  On the other hand, a remedy has been presented from the perspective of parallel computing devices, wherein such a distribution of a global signal proves costly.  Thus, the proposal is to allow cells to determine their own update time through local interactions (see \cite{fates2013} and the references therein, including \cite{nehaniv}).  Although the storage cost is higher, this scheme points towards a more natural form of computation.

In \cite{Wol1}, Wolfram explored synchronous CA on a one-dimensional lattice, where cells take the Boolean states 1 or 0.  The CA state of cell $i$ was dependent on the states of three connected cells, called the \emph{neighbourhood of $i$}.  These three-cell neighbourhood CA were termed ``elementary
cellular automata" (ECA).  

The ECA may be regarded as special cases of random Boolean networks.  The different types of these networks were first classified by Gershenson in 2002 \cite{Ger1}.  The types considered included asynchronous random Boolean networks (ARBNs), in which nodes are selected to be randomly updated at each time step, and deterministic ARBNs (DARBNs), where the node to be updated is selected deterministically. Gershenson talked of DARBNs as being
more advantageous because of their modelling capabilities, which are more straightforward
than ARBNs that rely on the stochasticity of asynchronous phenomena. Gershenson further proposed DARBNs as better representations of the famous genetic
regulatory networks of Kauffman \cite{Kauf1}, as they are asynchronous but do not rely on stochastic
methods.

Following on from Gershenson's idea of using determinism as a more `model-friendly' form for
asynchrony, a goal of this paper is to exploit this avenue by presenting a new, deterministic framework for asynchrony.  No matter how well it matches the real
system, we claim that the essence of many interesting and important asynchronous
processes is lost by using probabilistic updating schemes.

As a real example, consider Figure~\ref{fig:seashell} again, showing a seashell pattern; it is interesting to see that the same pattern (and many other such seashell patterns, as well as growth processes such as snowflakes \cite{snowflake}) may be replicated quite accurately by a CA model.  Traditionally, some difference between the two versions - real and CA - would be ascribed to a fault, random or otherwise, such that a better approximation may be obtained by adding stochastic asynchrony in the CA rule.  As alluded to by Cornforth et al in \cite{Cor2}, such asynchrony
has tended to simplify these dynamics into a
probability (or the like) of cells updating their states. Inspired by the ideas of `local clocks over global clocks' in \cite{fates2013} and \cite{nehaniv}, we will get into the heart
of the matter and study the pattern of exchanges that takes place locally, that is, before any data is transferred between cells.

\subsection{Network Description of the Cellular Automaton Lattice}
\begin{figure}[!hbp]
\centering
  \includegraphics[scale=0.4,trim= 50 240 20 470,clip]{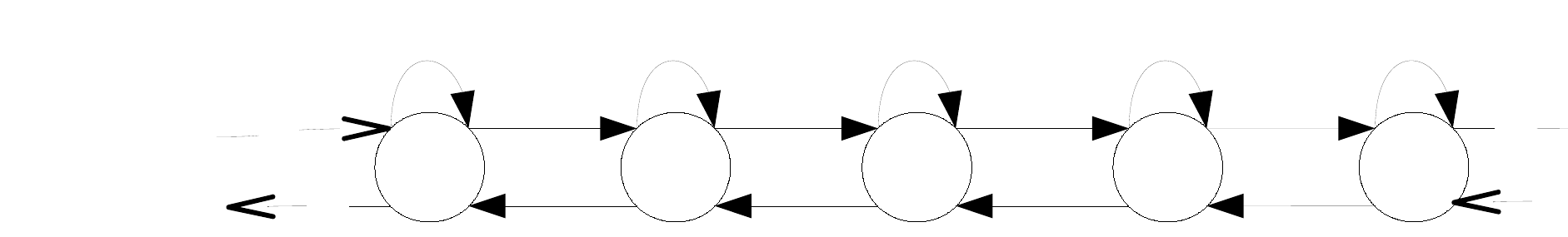} \\
  \caption{Digraph representing a regular 3-nbhd network}
  \label{fig:3nbhddigraph}
\end{figure}

Our work views each cell in a one-dimensional lattice as a processor
which receives input from its neighbourhood.  Having received this input, the processor computes its new state (as a function of the input states), then sends a corresponding output to its connected neighbours. This type of information exchange can be represented by the digraph in Figure~\ref{fig:3nbhddigraph}, where each node represents a processor and directed arcs between nodes indicate the direction of information transfer. 

Figure~\ref{fig:3nbhddigraph} shows three arcs pointing to each node, indicating that there are three
processors - therefore three neighbours (including $i$ itself) - sending information to each process $i$. Thus, the neighbourhood size $n$ of each node is $n = 3$; we also refer to such a neighbourhood as an \emph{$n$-neighbourhood} (or \emph{$n$-nbhd}). The figure particularly shows that each processor sends output to itself as well as to its left and right neighbours. We refer to this type of network as a \emph{regular $n$-nbhd network} or simply a regular network if $n$ is understood. Thus, the regular 3-nbhd network describes the lattice for the ECA, where cells are depicted by nodes.  For this reason, we use the terms \emph{cell} and \emph{node} to mean the same thing.

It is also assumed that there is a processing delay associated with
each nodal computation of CA state.  Moreover, we incorporate a transmission delay, which is the time taken for a CA state to be transmitted to other cells that require it.  These
two parameters are the means by which we obtain asynchrony: a divergence
from classical (synchronous) CA models since there the computations
are assumed to occur instantaneously.  

\subsection{Contents}
In Section~\ref{sec:syncCA}, we review the ECA, a classical synchronous model
of CA, which is later used to present our asynchronous model.  In Section~\ref{sec:async}, we show how max-plus algebra provides a natural way to mathematically model asynchronous CA.  By covering related graph theoretical techniques and known proven results, we also show that this asynchronous system is periodic, characterised by the connectivity of the underlying network.  It will be seen from the depth of theory covered that a max-plus algebraic model of asynchrony is more bespoke and addresses \emph{all} the intricacies of the internal dynamics within a cell (not just the external).  Nevertheless, those already familiar with max-plus algebra may skim over Sections \ref{sub:preliminaries} and \ref{sub:async}, noting the few places where we mention the link to our system.   The ``contour plot" is introduced as a framework for this asynchrony  in Section~\ref{sub:contour}.  This is followed by Section~\ref{sec:OutlineCA}, which shows the effect of the max-plus asynchrony on cellular automata.  We finish with concluding remarks.

\section{Synchronous Cellular Automata}\label{sec:syncCA}
Let $s_i$ denote the state of cell $i$. The index $i$ denotes the position of the cell in the one-dimensional Euclidean plane, so that the state of the CA at a given time $t\in\R$ is represented by the string $s_1(t)s_2(t)\cdots s_N(t)$ or the vector $\mathbf{s}(t) = (s_1(t),s_2(t),...,s_N(t))$, where $N$ is the size of the lattice.  The ECA assigns a symmetrical
neighbourhood of three nodes, where node $i$ is contained in its own
neighbourhood (as in Figure~\ref{fig:3nbhddigraph}).  Here, the \emph{CA rule} is a function
$f:\{0,1\}^3\rightarrow \{0,1\}$ given by $s_i(t+1)=f(s_{i-1}(t),s_i(t),s_{i+1}(t))$.  An example
of such a CA rule is the following.
\begin{equation} \label{equ:CArule150}
s_i(t+1)=\sum_{j=i-1}^{j=i+1}s_j(t)\quad \mod 2.
\end{equation}
The
rule in equation~(\ref{equ:CArule150}) is named ECA rule 150 by Wolfram
\cite{Wol2} and this is how we refer to it throughout the
report.

\begin{figure}[!hbt]
\centering

\includegraphics[scale=0.32, trim=310 45 330 20, clip]{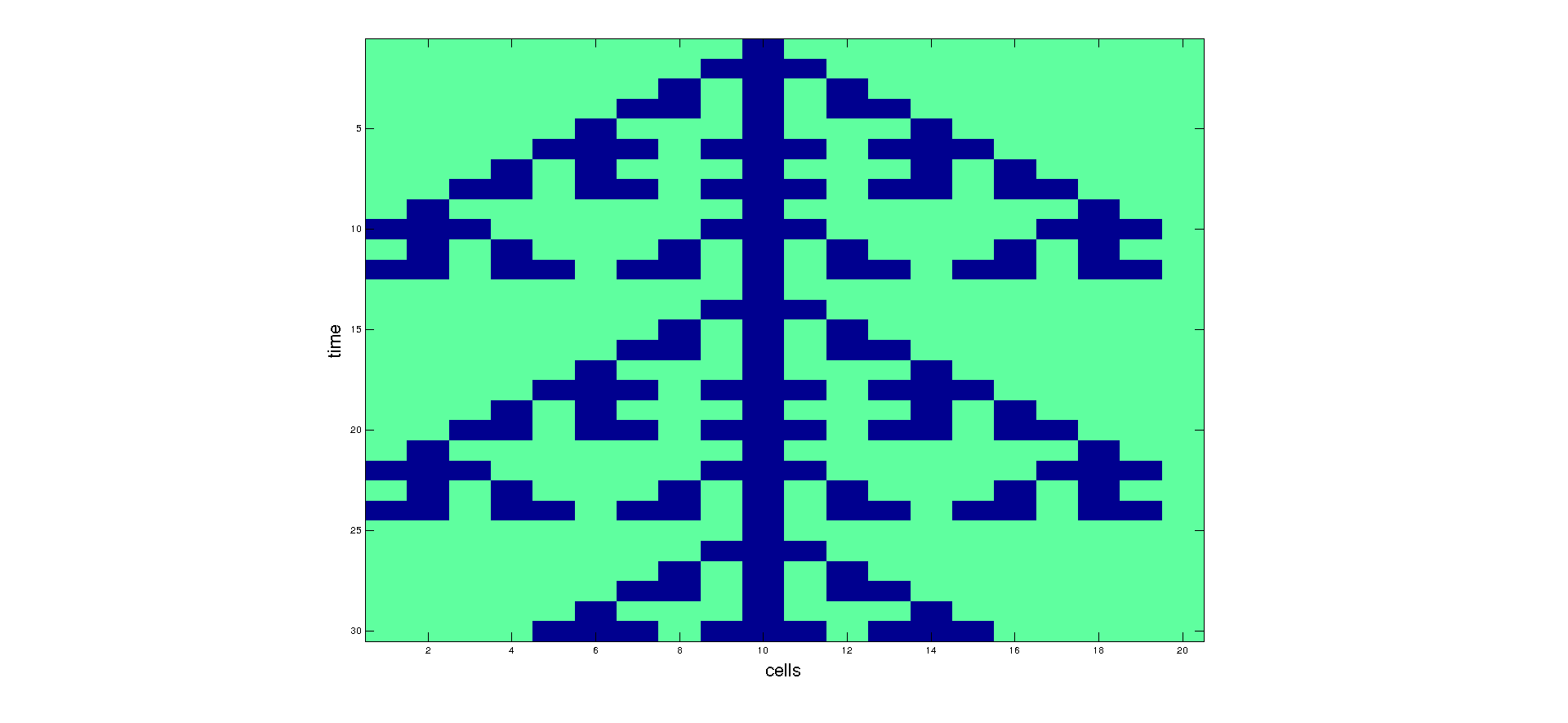}

\caption{CA pattern produced by ECA rule 150. Time travels
vertically down, cells are labelled on the $x$-axis.}\label{fig:carule150}
\end{figure}

Consider a regular 3-nbhd network of twenty cells (connected in a ring such that end
cells are adjacent).  If each cell is depicted by a square - coloured if
$s_i=1$, clear if $s_i=0$ - then the output produced after each of 30 iterations of equation~(\ref{equ:CArule150}) is shown in Figure~\ref{fig:carule150}.  For fixed $t$, the CA state $s_1(t)\cdots s_{20}(t)$ represents a horizontal line of cell states, and the initial CA state is
\begin{equation}
s_i(0)=\left\{ \begin{array}{ll}
                   1 & \text{if $i=10$}  \\
                   0 & \text{otherwise.}
                   \end{array} \right.
\end{equation}
Such an output as Figure~\ref{fig:carule150} is referred to as a \emph{space-time pattern} (or plot). It is evident that the absence of transmission and processing delays means that the $t^\text{th}$ update of each cell occurs at the same time as the $t^\text{th}$ update of every other cell, that is, synchronously.

\subsection{State Transition Graph}
Figure~\ref{fig:carule150} displays periodic behaviour. In other words, for each time step $t$, the CA yields a state $\ve{s}(t)$, which is
seen again after a few more time steps. Given the initial CA state
$\ve{s}(0)$ and the CA rule $f$, an \emph{orbit} of $\ve{s}(0)$ is
the sequence of states obtained by applying $f$ on $\ve{s}(0)$
sequentially.  If $f$ is applied $k$ times, we represent this as
$f^k(\ve{s}(0))=\underbrace{f(f(\cdots f}_{k\text{
times}}(\ve{s}(0))))$. We define periodic behaviour as follows.
\begin{defn}\label{def:CAperiod}
Consider the CA rule $f$ and network size $N$.  Let $\ve{s}(k)=f^k(\ve{s}(0))$ for all $k\geq0$, where $\ve{s}(k)$ is a CA state represented by a $1\times N$ vector.  For some $t\geq0$, if there exists a finite number $p\in\N$ such that $\ve{s}(t+p)=\ve{s}(t)$, then the set of states
\begin{displaymath}
\{\ve{s}(t),\ve{s}(t+1),\ldots,\ve{s}(t+p-1)\}
\end{displaymath}
is called \emph{a periodic CA orbit}, where $p$ is the \emph{CA period} of the orbit.
\end{defn}

\begin{figure}[!hbt]
\centering

 \includegraphics[scale=0.3, trim=0 265 0 165]{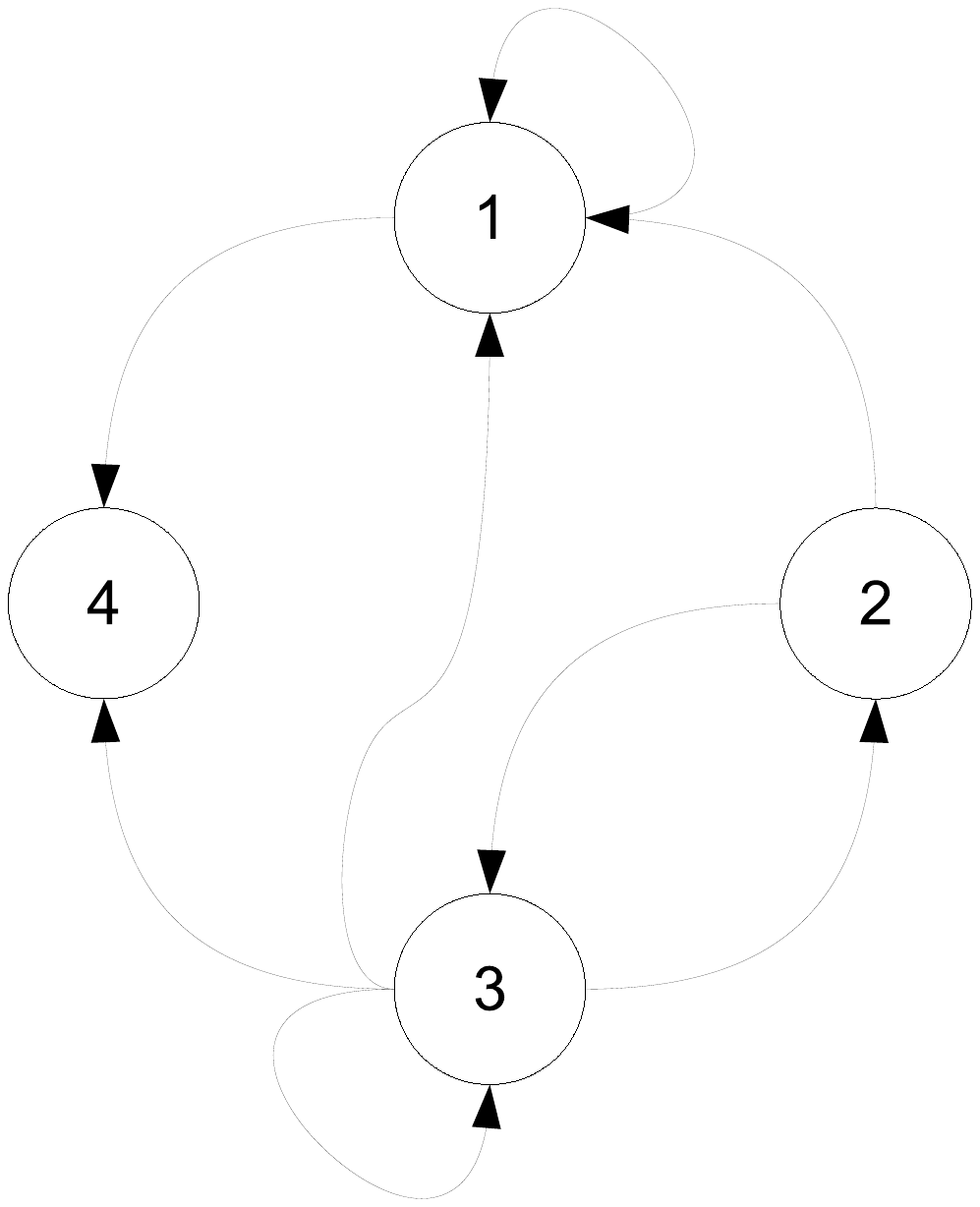}
 
  \caption{Size 4 network of CA cells.}
  \label{fig:STGexampledigraph}
\end{figure}

\begin{figure}[!hbt]
\centering

  \includegraphics[scale=0.5, trim=85 300 150 15,clip]{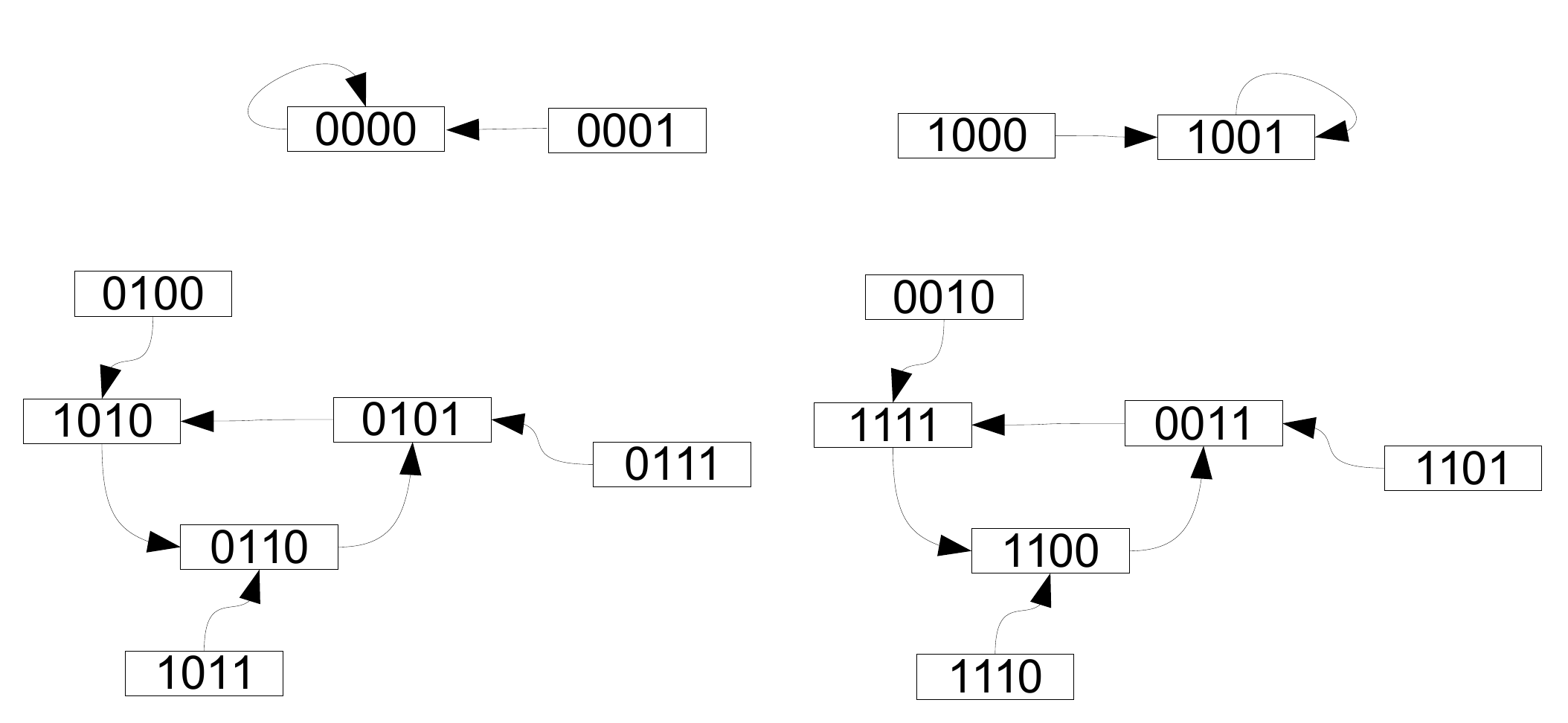}

  \caption{State transition graph of ECA rule 150 generalised to a lattice of size 4 (namely the network given in Figure~\ref{fig:STGexampledigraph}).}
  \label{fig:STGexample}
\end{figure}

As an example, consider the underlying network of size 4 as given in
Figure~\ref{fig:STGexampledigraph}.  We give the CA rule in words: the new state of each cell is the sum of the states ($\mod 2$) of its neighbourhood cells on the previous time step, where the neighbourhood of cell $i$ comprises those cells whose outgoing arc points to $i$.  (Note that this rule is an extension of ECA rule 150 to arbitrary lattices.)
For small $N$, as is the case here, it
is useful to represent each system state as a vertex in a digraph.  Thus,
there is an arc from CA state $\ve{s}^i$ to CA state $\ve{s}^{ii}$ if and
only if $f(\ve{s}^i)=\ve{s}^{ii}$.   The digraph is known as a \emph{state
transition graph} or STG for short.  For the system in question, the STG is given in
Figure~\ref{fig:STGexample}.  Each CA state is shown in string form,
where the $i^\text{th}$ digit represents the CA state of the $i^\text{th}$ node.

To determine the evolution of the CA, we can follow the arcs in Figure~\ref{fig:STGexample}.  It can be
seen that all initial CA states asymptotically evolve into four periodic orbits, represented as
circuits in the STG. Two of these (states 0000 and 1001) are
period-1 orbits (or \emph{fixed points} in
conventional dynamical systems language), and two are period-3
orbits; all other states are transient.  Such an STG is an artefact of a synchronous CA.

\section{Asynchronous Model}\label{sec:async}
\begin{figure}[!hbt]
\centering

  \includegraphics[scale = 0.26, trim= 230 155 190 130,clip]{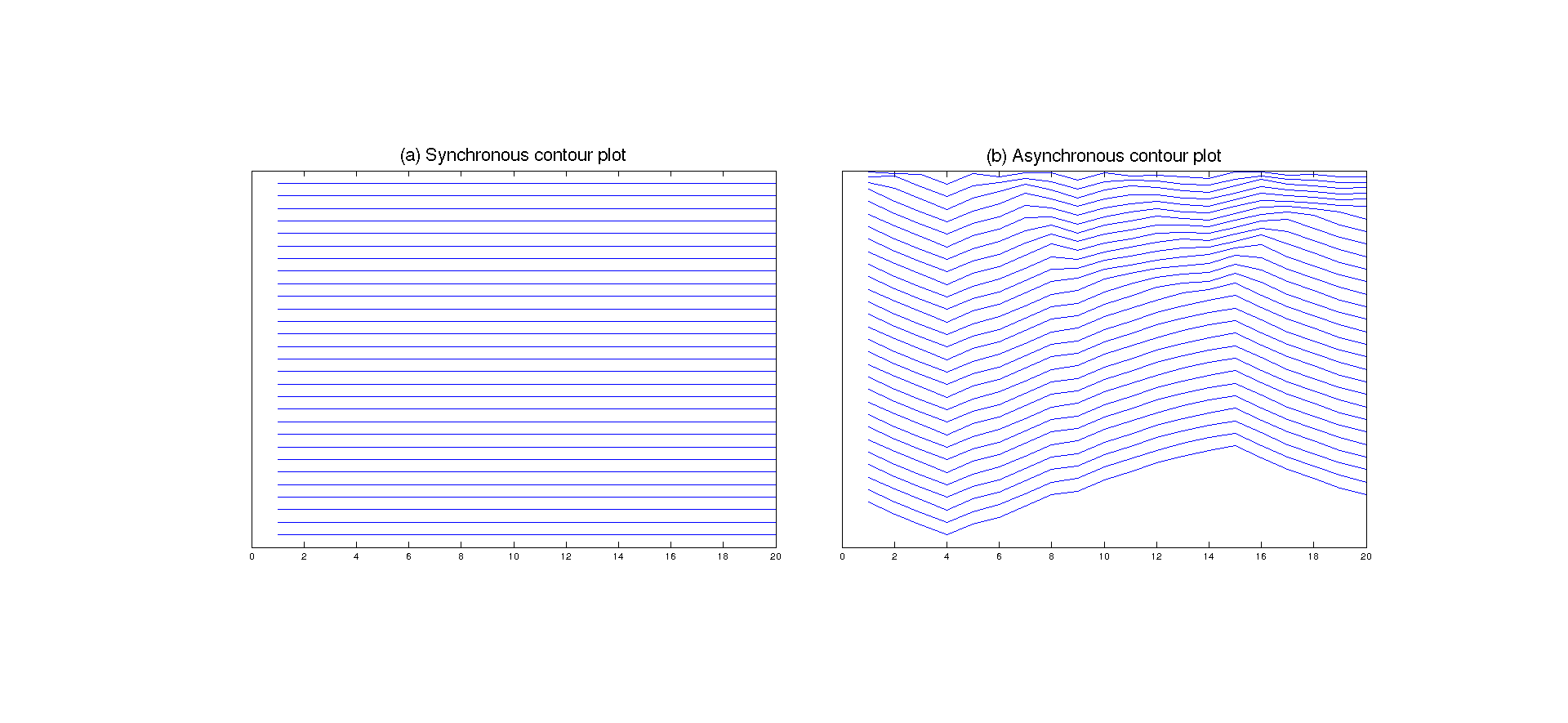} 

 \caption{(a) Synchronous and (b) asynchronous contour plots. The contours indicate update times of cells and act as a frame on which the CA may be evolved.  Both lattices are connected as a regular 3-nbhd network of 20 cells.  Time is on the vertical axis and the horizontal axis represents the cell positions.   In (a), contours are horizontal.  In (b), update $k$ of all cells is represented by contour $k$ (counting from the top).}
\label{fig:syncandasynccontourplots}
\end{figure}

Consider a synchronous cellular automaton.  Due to the synchrony, it is possible to draw horizontal lines in the corresponding space-time pattern such that each line represents the update times of all cells at some fixed time.
The CA is a discrete time dynamical system, so the horizontal
lines may be drawn in sequence, evenly spaced, as in
Figure~\ref{fig:syncandasynccontourplots}(a). We call such a space-time plot a \emph{contour plot}, and each horizontal line
is referred to as a \emph{contour}.  The contour plot may be thought
of as a frame on which the CA states are overlaid and simulated.
Now consider altering these contours so that cells do not
necessarily update synchronously.  The corresponding contour plot
may then look like Figure~\ref{fig:syncandasynccontourplots}(b), which shows
the contours having variable shapes - updates occur asynchronously.  We shall return to the contour plot after looking at how our asynchronous model may produce it.

\begin{figure}[!hbt]

\centering

  \includegraphics[scale=0.5,trim= 100 115 150 60,clip]{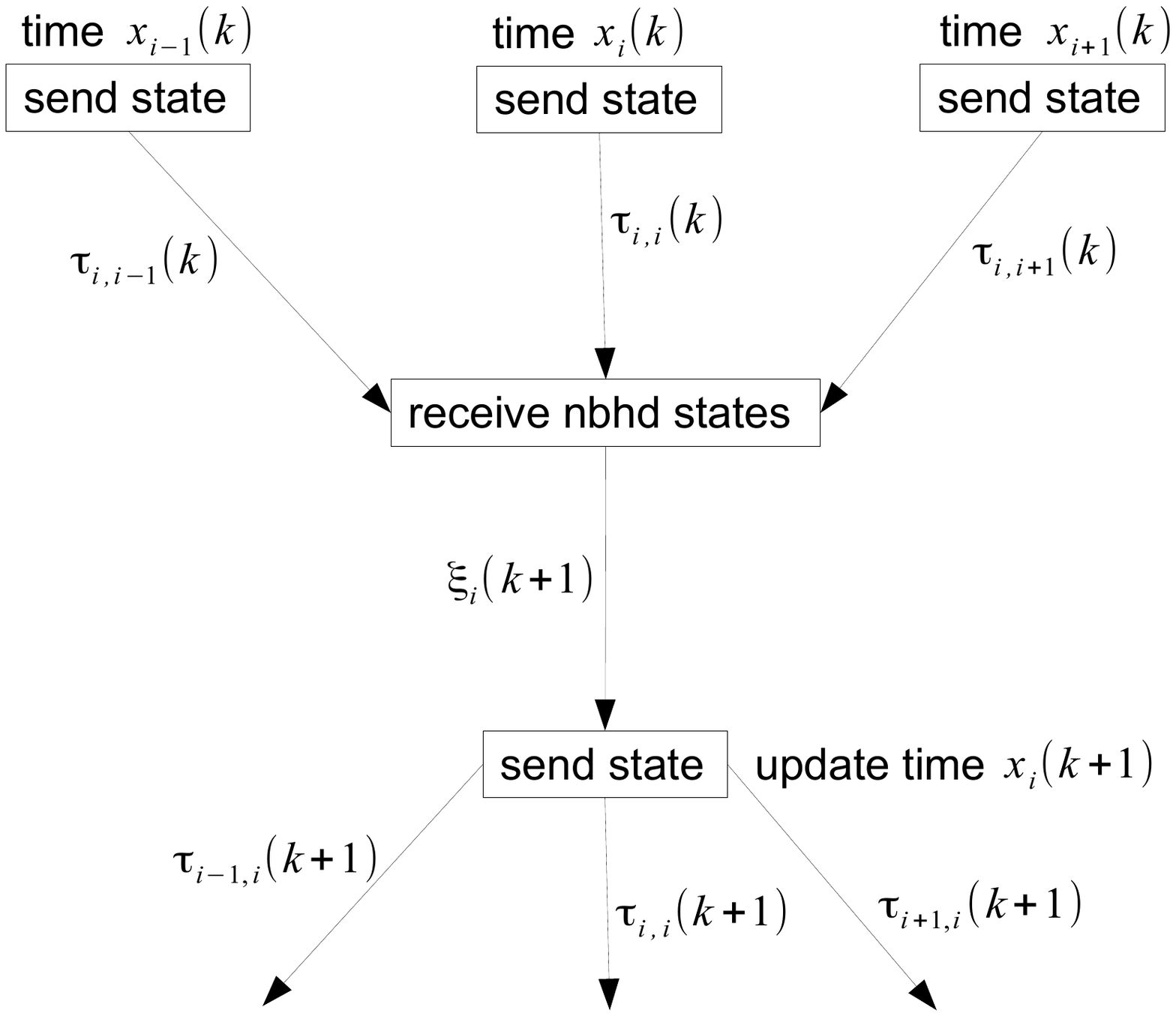} 

\caption{The processes internal to the $k^\text{th}$ state change at
cell $i$.  Real time travels vertically down. Arrows indicate the
destination of the sent state.  Labels on arrows indicate the time taken for the corresponding process.} \label{fig:statechange}
\end{figure}

We present the asynchronous model as follows.  Consider a cell $i$ contained in a regular $n$-nbhd network
of $N$ cells.  The cell carries a CA state (1 or 0) which changes with time depending on the rules that we
employ.  Thus, we can plot points on the real line corresponding to
when these changes occur.  The real line represents time and the
points are the \emph{update times} of the CA state.  Let $x_i(k)$
denote the $k^\text{th}$ update time for cell $i$.  We also refer to $k$ as a \emph{cycle number}.  Once each cell in the neighbourhood of cell $i$ has completed its $k^\text{th}$ cycle, it sends the updated state to $i$.  The transmission of such a state from cell $j$ to $i$ takes \emph{transmission time} $\tau_{ij}(k)$.  The update
of cell $i$ takes a \emph{processing time} and it is represented in the
$k^\text{th}$ cycle by $\xi_i(k)$.  If $n=3$, we have the
following iterative scheme for the $(k+1)^\text{th}$ update time of cell $i$.
\begin{eqnarray} 
x_i(k+1) =  \max  \{ & x_{i-1}(k)+\tau_{i,i-1}(k),x_i(k)+\tau_{i,i}(k),  \nonumber \\
          & x_{i+1}(k)+\tau_{i,i+1}(k)\}  +\xi_i(k+1) \label{equ:maxplusmodel1}
\end{eqnarray}
The above sequence of
interactions yielding a state change is depicted in
Figure~\ref{fig:statechange}.
Notice that we have now expanded on
the simpler notion of inter-cellular communication by
focussing closer on intra-cellular communication, that is, the communication within a cell itself.  This is the
key to our study of asynchrony, and it has arisen naturally from the requirement of a CA cell knowing the states of its neighbours.

We refer to these internal processes as \emph{events}.  In
Figure~\ref{fig:statechange}, there are two significant types of
events: ``receive" and ``send". The three times
$x_{i-1}(k)$, $x_i(k)$ and $x_{i+1}(k)$ are ``send" event
times, (i.e., when the corresponding CA states are sent).  The time
\begin{displaymath}
\max\{x_{i-1}(k)+\tau_{i,i-1}(k),x_i(k)+\tau_{i,i}(k),x_{i+1}(k)+\tau_{i,i+1}(k)\}
\end{displaymath}
is when node $i$ receives the aforementioned ``send" states; it is
therefore a ``receive" event.  Once received, node $i$ processes its
new CA state (by applying a CA rule on the received states); this takes time duration $\xi_i(k+1)$.  Once processed, node $i$ sends its
state to connected nodes at time $x_i(k+1)$; this is another ``send"
event. Note that we make no distinction between update times and events; both are the same.  We, however, formally refer to $x_i(k+1)$ as the update
time of $i$, as introduced earlier.  

\begin{note}
This idea of a network depicting events in a space-time diagram is not too dissimilar to Wolfram's \emph{causal network} in \cite[Chapter~9]{WolNKS}.  We shall delve into this relationship later.
\end{note}

\subsection{Preliminaries: Max-Plus Algebra and Graph Theory}\label{sub:preliminaries}
%
The $\max$ operation enables one to interpret equation~(\ref{equ:maxplusmodel1}) in \emph{max-plus algebra}.  This is useful because it converts a nonlinear system into a linear
system in this new algebra, which subsequently shares many important features with conventional linear algebra.  Particularly novel applications include the modelling of railway network timetables \cite{Heid}, manufacturing processes \cite{doust}, and even cellular protein production \cite{brackley}.  In presenting max-plus algebra, we borrow most notation and terminology from
\cite{Heid}.

Define $\varepsilon=-\infty$ and $e=0$, and denote by
$\mathbb{R}_\text{max}$ the set $\mathbb{R}\bigcup\{\varepsilon\}$.
For elements $a,b\in \mathbb{R}_\text{max}$, define operations
$\oplus$ and $\otimes$ by
\begin{displaymath}
a\oplus b = \max(a,b)\qquad \text{and}\qquad a\otimes b = a+b.
\end{displaymath}
The set $\mathbb{R}_\text{max}$ together with the operations
$\oplus$ and $\otimes$ is what we refer to as \emph{max-plus
algebra} and it is denoted by $\mathcal{R}_\text{max}=(\mathbb{R}_\text{max},\oplus,\otimes,\varepsilon,e)$.  $\varepsilon=-\infty$ is the `zero' (i.e., $\forall x\in\R_{\max}$, $\varepsilon\otimes
x=x\otimes \varepsilon=\varepsilon$ and
$\varepsilon\oplus x=x\oplus \varepsilon=x$), whilst $e=0$ is the `unit' element (i.e., $\forall x\in\R_{\max}$, $e\otimes x=x\otimes e=x$).

$\mathcal{R}_{\max}$ is associative and commutative over both operations
$\oplus$ and $\otimes$ while $\otimes$ is distributive over $\oplus$.  (In addition, $\oplus$ is idempotent in $\mathcal{R}_{\max}$, so that max-plus algebra is a commutative and idempotent semiring.)  

In this paper, we remove the dependence on $k$ of the processing and transmission times so that equation~(\ref{equ:maxplusmodel1}) is written 
\begin{equation}
x_i(k+1)=\max\{x_{i-1}(k)+\tau_{i,i-1},x_i(k)+\tau_{i,i},x_{i+1}(k)+\tau_{i,i+1}\}+\xi_i.
\end{equation}
We can now write this in max-plus algebra:
\begin{eqnarray}
x_i(k+1) =  \{ & (\tau_{i,i-1}\otimes x_{i-1}(k))\oplus(\tau_{i,i}\otimes x_i(k))\nonumber \\
          & \oplus(\tau_{i,i+1}\otimes x_{i+1}(k))\}  \otimes\xi_i \label{equ:maxplusmodel3} 
\end{eqnarray}
We often save space and clarify the presentation by omitting $\otimes$, much as in conventional algebra.  Thus, $x\otimes
y\equiv xy$ and equation~(\ref{equ:maxplusmodel3}) can be
written as
\begin{equation}
x_i(k+1)=\xi_i\{\tau_{i,i-1}x_{i-1}(k)\oplus\tau_{i,i}x_i(k)\oplus\tau_{i,i+1}x_{i+1}(k)\}
\end{equation}
Since $\otimes$ is distributive over $\oplus$, we can write this as
\begin{equation} \label{equ:maxplusmodel2}
x_i(k+1)=\xi_i\tau_{i,i-1}x_{i-1}(k)\oplus\xi_i\tau_{i,i}x_i(k)\oplus\xi_i\tau_{i,i+1}x_{i+1}(k)
\end{equation}
To represent a max-plus power, we follow from associativity of
$\otimes$ and define, for $x\in\R_{\max}$,
\begin{equation}
x^{\otimes n}\stackrel{\text{def}}{=}\underbrace{x\otimes x\otimes \cdots \otimes x}_\text{$n$ times}
\end{equation}
for all $n\in\N$ with $n\neq0$. For
$n=0$, we define $x^{\otimes0}=e$ ($=0$).  

Max-plus algebra extends naturally to matrices, and this allows
the concurrent modelling of the update times for all nodes. Denote
the set of $n\times m$ matrices with underlying max-plus algebra by
$\R_{\max}^{n\times m}$.  The sum of matrices $A,B \in\R_{\max}^{n\times m}$, denoted by
$A\oplus B$, is defined by
\begin{equation}
[A\oplus B]_{ij} = a_{ij}\oplus b_{ij}
\end{equation}
where $a_{ij}=[A]_{ij}$ and $b_{ij}=[B]_{ij}$.  In the same vein, for matrices
$A\in\R_{\max}^{n\times l}$ and $B\in\R_{\max}^{l\times m}$, the
matrix product $A\otimes B$ is defined by
\begin{equation}
[A\otimes B]_{ij} = \bigoplus_{k=1}^l a_{ik}\otimes b_{kj} = \max_{k\in\{1,\ldots,l\}} \{a_{ik}+b_{kj}\}
\end{equation}
For $\alpha\in\R$, the scalar multiple $\alpha\otimes A$ is defined
by
\begin{equation}
[\alpha\otimes A]_{ij}=\alpha\otimes a_{ij}.
\end{equation}


As in classical matrix manipulation, the max-plus matrix addition
$\oplus$ is associative and commutative, whilst the matrix product
$\otimes$ is associative and distributive with respect to $\oplus$;
it is usually not commutative.   Similarly,
the operation $\otimes$ has priority over $\oplus$.

The elements of $\R_{\max}^n
\stackrel{\text{def}}{=}\R_{\max}^{n\times 1}$ are called
\emph{vectors}.  A vector is usually written in bold, as in $\ve{x}$, whilst the vector with all elements equal to $e$ is called
the \emph{unit vector} and is denoted by $\ve{u}$.

We are now able to define matrix-vector products.  The product $A\otimes
\ve{x}$, where $A\in\R_{\max}^{n\times m}$ and
$\ve{x}\in\R_{\max}^n$ is
\begin{equation}
[A\otimes \ve{x}]_i = \bigoplus_{k=1}^m a_{ik}\otimes x_k  = \max_{k\in\{1,\ldots,m\}} \{a_{ik}+x_k\}.
\end{equation}
Moreover, for the square matrix $A\in\R_{\max}^{n\times n}$, denote
the $k^\text{th}$ power of $A$ by $A^{\otimes k}$, defined by
\begin{equation}\label{equ:matrixpower}
A^{\otimes k} \stackrel{\text{def}}{=} \underbrace{A\otimes A\otimes \cdots \otimes A}_\text{$n$ times}
\end{equation}
for all $k\in\N$ with $k\neq0$. For $k=0$, we set
$A^{\otimes0}\stackrel{\text{def}}{=}E(n,n)$, the identity matrix whose diagonal elements equal $e$ and all of whose other elements are $\varepsilon$.

Having established the preliminaries above, a system of $N$ such
equations as (\ref{equ:maxplusmodel2}) can now be given in the form
\begin{equation} \label{equ:xk+1P}
\ve{x}(k+1)=P\otimes \ve{x}(k)
\end{equation}
where $\ve{x}(k)=(x_1(k), x_2(k),\ldots,x_N(k))^\top$.  $P$ is
the $N\times N$ matrix defined by $A_\xi\otimes T$, where
\begin{displaymath}
A_\xi = \left(
  \begin{array}{cccc}
    \xi_1 & \varepsilon & \cdots & \varepsilon \\
    \varepsilon & \xi_2 & \cdots & \varepsilon \\
    \vdots &  & \ddots & \vdots \\
    \varepsilon & \varepsilon & \cdots & \xi_N \\
  \end{array}
\right)
\end{displaymath}
and
\begin{displaymath}
T =\left (\begin{array}{ccccccc}
               \tau_{11} & \tau_{12} & \varepsilon & \varepsilon & \cdots & \varepsilon & \tau_{1N}\\
               \tau_{21} & \tau_{22} & \tau_{23} & \varepsilon & \cdots & \varepsilon & \varepsilon \\
               \vdots &  &  & \ddots & & & \vdots \\
               \varepsilon & \varepsilon & \cdots & \varepsilon & \tau_{N-1,N-2} & \tau_{N-1,N-1} & \tau_{N-1,N} \\
               \tau_{N,1} & \varepsilon & \cdots & \varepsilon & \varepsilon & \tau_{N,N-1} & \tau_{NN}
             \end{array}
 \right ).
\end{displaymath}
$A_\xi$ is referred to as the \emph{processing matrix} and $T$ is the \emph{transmission matrix}. We call equation~(\ref{equ:xk+1P}) a \emph{max-plus system} (of dimension $N$) where
the vector $\ve{x}(k)$ is the \emph{state of the system}.  $P$ is called the \emph{timing dependency matrix}\footnote{This name (along with the later \emph{timing dependency graph}) is inspired by its use in another novel application of $\mathcal{R}_{\max}$ to the timing of digital hardware in \cite{broomheaddigital}.} of the network of cells.  NB:  Here, the term ``state" refers to update time and is not to be confused with ``CA state".  Nevertheless, the context should make this distinction clear.

To a network of cells, we associate a
\emph{digraph}, such as in Figure~\ref{fig:3nbhddigraph}.  In general, we define a digraph as $\mathcal{G}=(V,E)$, consisting of a set $V$ and a
set $E$ of ordered pairs $(a,b)$ of $V$.  (Often, we refer to the digraph simply as a graph which, in turn, also refers to the ``network" of our application.)  The elements of $V$ are
called \emph{vertices} or \emph{nodes} and those of $E$ are \emph{arcs}. 

An arc $(a,b)$ is also denoted $ab$, and we refer to an
arc $aa$ as a \emph{self-loop}.  For the arc $ab$, $a$ is the
\emph{start node} and $b$ is the \emph{end node};  $a$ is also referred to as the \emph{predecessor} node of $b$ whilst $b$ is the \emph{successor} of $a$. 

By assigning real numbers (called \emph{weights}) to the arcs of a graph $\mathcal{G}=(V,E)$, we obtain a
\emph{weighted graph}.   The \emph{weighted adjacency
matrix of $\mathcal{G}$ over $\R_{\max}$} is the matrix $W\in \R_{\max}^{n\times
n}$ whose $(i,j)^\text{th}$ entry $w_{ij}$ is non-zero ($\neq \varepsilon$) if and only if $j$ is a predecessor of $i$.  We also refer to $W$ as a \emph{max-plus adjacency matrix}.  Given $W\in\R_{\max}^{n\times
n}$, we denote the associated network as $\mathcal{G}(W)$.  
Thus, $P$ in equation~(\ref{equ:xk+1P}) is a max-plus adjacency matrix of the regular $3$-nbhd network since the network is exactly that shown in Figure~\ref{fig:3nbhddigraph} augmented with arcweights $\xi_i\tau_{ij}$.  In fact, since our particular max-plus system concerns update times, we call $\mathcal{G}(P)$ the \emph{timing dependency graph} of the system.  The neighbourhood of a node may now be defined in terms of the adjacency matrix.
\begin{defn} \label{def:nbhd}
Let $W$ be the max-plus adjacency matrix of a digraph of connected cells. The
\emph{neighbourhood} of $i$ is $\mathcal{N}_i=\{j|w_{ij}\neq\varepsilon\}$.
\end{defn}

\begin{defn}
Let $p=\{a_1,a_2,\ldots,a_n\}$ be a sequence of arcs.  If there are vertices $v_0,v_1,\ldots,v_n$ (not necessarily
distinct) such that $a_j=v_{j-1}v_j$ for
$j=1,\ldots,n$ then $p$ is called a \emph{walk} from $v_0$ to
$v_n$.  A walk for which the $a_j$ are distinct is called a \emph{path}.  Such a path is said to consist of the nodes $v_0,v_1,\ldots,v_n$ and
to have \emph{length} $n$, which is denoted $|p|_l=n$.

If $v_n=v_0$, then the path is called a \emph{circuit}.  If the
nodes in the circuit are all distinct (i.e., $v_i\neq v_k$ for $i\neq
k$), then it is called an \emph{elementary} circuit.
\end{defn}
We define the
\emph{weight $|p|_w$ of a path $p$} as the sum of the weights of all arcs
constituting the path.  
The \emph{average weight} of $p$ is $\frac{|p|_w}{|p|_l}$.
For a circuit, we refer to this quantity as the \emph{average
circuit weight}.
\begin{defn}
For a graph $\mathcal{G}=(V,E)$,
node $j\in V$ is said to be \emph{reachable} from node $i\in V$, denoted $i\rightarrow j$, if there
exists a path from $i$ to $j$.  Graph $\mathcal{G}$ is
\emph{strongly connected} if $i\rightarrow j$ for any two nodes $i,j\in V$.
\end{defn}
Moreover, matrix $A\in\R_{\max}^{n\times n}$ is called
\emph{irreducible} if $\mathcal{G}(A)$ is strongly connected; if a
matrix is not irreducible, it is called \emph{reducible}.  Thus, in our max-plus system of equation~(\ref{equ:xk+1P}), matrix $P$ is irreducible.

\begin{defn}
Denote the \emph{cyclicity} of a graph $\mathcal{G}$ by $\sigma_{\mathcal{G}}$.
\begin{itemize}
\item If $\mathcal{G}$ is strongly connected, then $\sigma_{\mathcal{G}}$ equals the greatest common divisor of the lengths of all elementary circuits in $\mathcal{G}$.  If $\mathcal{G}$ consists of only one node without a self-loop, then $\sigma_{\mathcal{G}}$ is defined to be one
\item If $\mathcal{G}$ is not strongly connected, then $\sigma_{\mathcal{G}}$ equals the least common multiple of the cyclicities of all maximal strongly connected subgraphs (MSCSs) of $\mathcal{G}$. (See Appendix~\ref{app:MSCS} for definition of MSCSs.)
\end{itemize}
\end{defn}
\begin{defn} \label{def:cyclicity}
Let $A\in\R_{\max}^{n\times n}$ be irreducible.  The \emph{cyclicity
of $A$}, denoted $\sigma(A)$, is defined as the cyclicity of the
critical graph of $A$.
\end{defn}
When the matrix is understood, the cyclicity is also denoted by $\sigma$.


\subsection{Asymptotic Behaviour of the Max-Plus System}\label{sub:async}
Let $\ve{x}(0)$ represent the initial state of all cells.  Then we
can rewrite equation~(\ref{equ:xk+1P}) as
\begin{equation}
\ve{x}(k+1) = P\otimes P\otimes \cdots P \otimes \ve{x}(0) =P^{\otimes(k+1)}\otimes\ve{x}(0)
\end{equation}
or, equivalently,
\begin{equation}\label{equ:xkPkx0}
\ve{x}(k)=P^{\otimes k}\otimes \ve{x}(0).
\end{equation}
Given $\ve{x}(0)$, the sequence of vectors
$\{\ve{x}(k):k\in\N_0\}$, obtained by iterating equation~(\ref{equ:xkPkx0}), is
referred to as the \emph{orbit} of $\ve{x}(0)$.  A study of such
sequences is provided in greater scope in \cite[Chapters~3 and 4]{Heid}.  We detail the topics relevant for this work, with some well-known results expanded upon in Appendix~\ref{app:asymptotic}.

\begin{defn}
Let $A\in\R_{\max}^{n\times n}$.  If $\lambda\in\R_{\max}$ is a
scalar and $\ve{v}\in\R_{\max}^n$ is a vector that contains at least
one finite element such that
\begin{equation}
A\otimes \ve{v}=\lambda\otimes\ve{v},
\end{equation}
then $\lambda$ is called
an \emph{eigenvalue} of $A$ and $\ve{v}$ is an \emph{eigenvector} of
$A$ associated with eigenvalue $\lambda$.
\end{defn}
 
For a system with irreducible $P$ such as equation~(\ref{equ:xk+1P}) (or~(\ref{equ:xkPkx0})), it turns out that there
is only one eigenvalue, and it is equal to the maximal average weight of elementary circuits in $\mathcal{G}(P)$ (See Appendix~\ref{app:asymptotic}, Theorem~\ref{thm:evalueirr}).  Such circuits with maximal average weight are called \emph{critical} and the \emph{critical
graph} of $P$ is the graph
consisting only of critical circuits
in $\mathcal{G}(P)$.  

\begin{defn}\label{def:periodic}
Let $A\in\R_{\max}^{n\times n}$.  For some $k\geq0$, consider the
set of vectors
\begin{displaymath}
\ve{x}(k),\ve{x}(k+1),\ve{x}(k+2),\ldots \in\R_{\max}^n
\end{displaymath}
where $\ve{x}(k)=A^{\otimes k}\otimes\ve{x}(0)$ for all $k\geq0$. The set
is called a \emph{(periodic) regime} if there exists $\mu\in\R_{\max}$
and a finite number $\rho\in\N$ such that
\begin{displaymath}
\ve{x}(k+\rho)=\mu\otimes\ve{x}(k).
\end{displaymath}
The \emph{period} of the regime is $\rho$.
\end{defn}
Suppose the initial vector $\ve{x}(0)$ is an eigenvector of $A$.  Then $\ve{x}(k+1)=\lambda\otimes\ve{x}(k)$ for $k\geq0$, so that the period is one.  Thus, larger periods are obtained when the system is not initialised to an eigenvector. The
remainder of Section~\ref{sec:async} explores this in more detail.  We start with the following crucial theorem of max-plus
algebra.
\begin{theorem}\label{thm:Alambdacyc}
Let $A\in\R_{\max}^{n\times n}$ be an irreducible matrix with
eigenvalue $\lambda$ and cyclicity $\sigma$.  Then there is a $k_\star$
such that
\begin{displaymath}
A^{\otimes(k+\sigma)}=\lambda^{\otimes\sigma}\otimes A^{\otimes k}
\end{displaymath}
for all $k\geq k_\star$.
\begin{proof}
See \cite[Theorem 3.9]{Heid}.
\end{proof}
\end{theorem}
Let $\sigma=\sigma(P)$.  For an indication of the asymptotic
behaviour of our system, we apply this theorem to observe the
state at epoch $k+\sigma$ for $k\geq k_\star$:
\begin{eqnarray}
\ve{x}(k+\sigma) &=& P^{\otimes (k+\sigma)} \otimes\ve{x}(0) \nonumber \\
                        &=& \lambda^{\otimes\sigma}\otimes P^{\otimes k} \otimes \ve{x}(0) \nonumber \\
                        &=& \lambda^{\otimes\sigma} \otimes \ve{x}(k)
\end{eqnarray}
where $\lambda^{\otimes\sigma}$ is read as $\lambda\times\sigma$ in
terms of classical algebra.  This guarantees the periodic behaviour of the max-plus system, where $\sigma$ is the upper bound on the period.  More specifically, the period $\rho$ is dependent on the choice of $\ve{x}(0)$ and, since we have seen that the system must also be periodic with period $\sigma$, we have that $\rho$ is a factor of $\sigma$.  We saw earlier that if $\ve{x}(0)$ is not an eigenvector, then $\rho>1$; we now also know that such a period $\rho$ will not be larger than $\sigma$.

We can use the above to show that the vectors $\ve{x}(k)$ in a regime turn out to be eigenvectors of
$P^{\otimes\sigma}$ associated with eigenvalue
$\lambda\times\sigma$:
\begin{eqnarray}
\lambda^{\otimes\sigma} \otimes \ve{x}(k) &\stackrel{\text{from above}}{=}& P^{\otimes (k+\sigma)} \otimes\ve{x}(0) \nonumber \\
 &=& P^{\otimes\sigma}\otimes P^{\otimes k}\otimes\ve{x}(0) \nonumber \\
 &=& P^{\otimes\sigma}\otimes \ve{x}(k).
\end{eqnarray}
In fact, given $\ve{x}(0)$ and corresponding period $\rho$, vectors in a regime are also eigenvectors of $P^{\otimes\rho}$; this can be shown in the same way as above. 

We now define a measure for the average delay between consecutive event times $x_i(k)$ and $x_i(k+1)$.
\begin{defn}\label{def:cycletimeasymp}
Let $\{x_i(k):k\in\N\}$ be an orbit of $x_i(0)$ in $\R_{\max}$.
Assuming that it exists, the quantity $\chi_i$, defined by
\begin{displaymath}
\chi_i = \lim_{k\rightarrow\infty}\frac{x_i(k)}{k}
\end{displaymath}
is called the \emph{cycletime} of $i$.
\end{defn}  

For an irreducible system such as ours, the vector $\ve{\chi}=(x_1,x_2,\ldots,x_N)$ of cycletimes is unique (See Appendix~\ref{app:asymptotic}, Theorem~\ref{thm:cycletimeindependent}).  Moreover, the irreducibility of $P$ ensures that each element of $\ve{\chi}$ is the same - specifically the eigenvalue of $P$, that is, 
\begin{displaymath}
\lim_{k\rightarrow\infty}\frac{x_i(k)}{k}=\lambda
\end{displaymath}
for any initial condition $\ve{x}(0)\in\R^n$. (See Appendix~\ref{app:asymptotic}, Lemma~\ref{lem:allnodessamecycle}.)

Thus, since $\ve{\chi}$ is independent of the initial condition, we relate it to our timing dependency graph by calling it the \emph{cycletime vector of $P$}.  $P$ is irreducible in this paper, so we let $\ve{\chi}=\chi_i$ for any $i$, and therefore refer to the cycletime vector of irreducible $P$ simply as the cycletime of $P$. 

As a compact summary of this subsection, we have shown that asynchrony due to an irreducible max-plus system always leads to periodic behaviour, and it is characterised by the circuit(s) in $\mathcal{G}(P)$ with largest average weight.  Section~\ref{sub:contour} addresses the impact of the above theory on the asynchronous time framework of the contour plot introduced earlier.

\subsection{The Contour Plot}\label{sub:contour}
Figure~\ref{fig:statechange} can be seen as a \emph{Hasse diagram} of events.  We say that events are \emph{causally related} if they are contained in the same \emph{chain}.  For example, the
``send" events at times $x_{i-1}(k-1)$, $x_i(k-1)$ and
$x_{i+1}(k-1)$ are not causally related since they are not contained in a chain - there is no path of directed arcs
connecting any of the three events.  Consequently, these three events form an \emph{antichain}.  (For a formal definition of a Hasse diagram, including chains and antichains, we refer the reader to Appendix~\ref{app:hasse}.)  By connecting those
elements in the same antichain, we obtain a piecewise linear
plot of the vector $\ve{x}(k)$, which we define next.
\begin{defn}\label{def:contour}
Consider the vector $\ve{x}(k)$.  A \emph{contour} is the plot
obtained by connecting $(i,x_i(k))$ to $(i+1,x_{i+1}(k))$ with a
straight line for each $i$, ($i=1,\ldots,N$).
\end{defn}
Creating a
contour for each $k$ gives a pictorial representation of vectors
$\ve{x}(k)$ as a function of $k$. We call this a \emph{contour plot}.
Figure~\ref{fig:syncandasynccontourplots} displays the contour plots of a size 20 system, where the
sequence $\{\ve{x}(0),\ve{x}(1),\ve{x}(2),\ldots\}$ represents
the contours (counting $k$ from the top).  For this reason, we interchangeably refer to the vector $\ve{x}(k)$ by ``the $k^\text{th}$
contour" from now on.  Between successive contours, we can imagine there being drawn the
internal processes of those in Figure~\ref{fig:statechange}.

Consider an example system with $\ve{x}(0)=\ve{u}$ that yields the following periodic behaviour: $\ve{x}(k+1)=5\otimes\ve{x}(k)$ for $k\geq 3$.  The contours for this system would represent vectors in the periodic regime $\{\ve{x}(k)|\ve{x}(k+1)=5\otimes\ve{x}(k),\ve{x}(0)=\ve{u},k\geq 3\}$.  The period of a regime and cyclicity are related by
$1\leq\rho\leq\sigma$. Thus, if $\sigma=1$ in this
example, then $\rho=\sigma=1$, so that no other period can be obtained
for all initial states $\ve{x}(0)$.  Therefore, each
contour in the contour plot has the same shape (separated by 5 time units) as $k\rightarrow\infty$; we call this a \emph{limiting shape} of the contours or a \emph{limiting contour}.  For a larger
period, we obtain a different set of limiting contours.  In particular, we obtain different limiting contour plots for the cases
$\rho=1,2,\ldots,\sigma$, each dependent on the choice of
$\ve{x}(0)$.  We will formalise this in Section~\ref{sec:OutlineCA}.

The idea of a limiting shape in contours suggests a
change of coordinates:  Given the irreducible matrix
$P\in\R_{\max}^{N\times N}$ with eigenvalue $\lambda$, let
\begin{equation}\label{equ:movingframe}
\ve{x}(k)=\lambda^{\otimes k}\otimes\ve{y}(k).
\end{equation}
We can think of $\lambda^{\otimes k}$ as a diagonal matrix or the
product of $\lambda^{\otimes k}$ and the identity matrix $E(N,N)$.  The
advantage of this is that such a diagonal matrix is invertible, its
inverse being the diagonal matrix with diagonal entries equal to
$\lambda^{\otimes-k}$.  Using this property, we rearrange equation
(\ref{equ:movingframe}) to obtain
\begin{equation}
\ve{y}(k)=\lambda^{\otimes -k}\otimes\ve{x}(k).
\end{equation}
In other words, $\ve{y}(k)$ is the limit to which the vectors
$\ve{x}(k)-\lambda k$ tend to as $k\rightarrow \infty$.  By studying
the asymptotic behaviour of $\ve{y}(k)$ itself, we can deduce the shape of the limiting contour.  

The original system follows the recurrence relation
$\ve{x}(k+1)=P\otimes\ve{x}(k)$ for some $\ve{x}(0)\in\R_{\max}^N$.
Substitute equation~(\ref{equ:movingframe}) into this to obtain
\begin{equation}
\lambda^{\otimes(k+1)}\otimes\ve{y}(k+1)=P\otimes\lambda^{\otimes
k}\otimes\ve{y}(k).
\end{equation}
Interpreting $\lambda^{\otimes k}$ as a diagonal matrix again yields
\begin{eqnarray}
\ve{y}(k+1) &=& \lambda^{\otimes -(k+1)}\otimes P\otimes\lambda^{\otimes k}\otimes\ve{y}(k) \nonumber \\
                &=& \lambda^{\otimes -1}\otimes P\otimes\ve{y}(k) \nonumber \\
                &=& \hat{P}\otimes\ve{y}(k) \label{equ:yPhat}
\end{eqnarray}
where $\hat{P}=\lambda^{\otimes -1}\otimes P$ represents the
\emph{normalised matrix} of $P$, equivalently obtained by
subtracting the eigenvalue of $P$ from each of its entries.  The
communication graph of $\hat{P}$ is the same as that for $P$ (but with different arcweights) so that $\hat{P}$ is also irreducible.  However, the maximum
average circuit weight of $\mathcal{G}(\hat{P})$, hence the
eigenvalue of $P$, is zero.

Moreover, it can be shown that $P$ and $\hat{P}$ have the same cyclicity (i.e.,
$\sigma(\hat{P})=\sigma(P)=\sigma$).  Theorem \ref{thm:Alambdacyc}
tells us of the asymptotic behaviour of the powers of an irreducible
matrix. Apply this to $\hat{P}$ to obtain $\hat{P}^{\otimes (k+\sigma)} =0^{\otimes\sigma}\otimes\hat{P}^{\otimes k} =\hat{P}^{\otimes k}$ for $k\geq k_\star$.  Thus, using equation~(\ref{equ:yPhat}),
\begin{eqnarray}
\ve{y}(k+\sigma) &=& \hat{P}^{\otimes (k+\sigma)}\otimes\ve{y}(0) \nonumber \\
                     &=& \hat{P}^{\otimes k} \otimes\ve{y}(0) \nonumber \\
                     &=& \ve{y}(k).
\end{eqnarray}
So the limiting contour $\ve{y}(k)$ is periodic with period
$\sigma(P)$.  In fact, like $\ve{x}(k)$, $\ve{y}(k)$ is periodic with period $\rho$, dependent on $\ve{x}(0)$, where $\rho$ is a factor of $\sigma$. Note that this period now conforms with the
traditional dynamical systems definition of a period in that the
sequence $\{\ve{y}(k)|k\in\N_0\}$ is not monotonically increasing, in contrast to
the original sequence $\{\ve{x}(k)|k\in\N_0\}$.  

Loosely speaking, Section~\ref{sub:contour} has shown that there is no unique shape to the limit of a contour
plot.  This yields an interesting feature of the max-plus asynchronous model: the asynchrony is related not only to the timing dependency graph but also to the system's starting point in time.

\subsection{The Eigenspace in Max-Plus Algebra} 

The set of all eigenvectors of $A\in\R_{\max}^{n\times n}$ associated to eigenvalue
$\lambda$ is the \emph{eigenspace} of $A$. For max-plus asynchrony, the importance of the eigenspace is demonstrated through its links to the contour plot, as shown next.

Consider the recurrence
relation $\ve{x}(k+1)=P\otimes\ve{x}(k)$ and its corresponding contour plot.  In Section~\ref{sec:async}, we
established that the vectors in a regime are eigenvectors of $P^{\otimes\sigma(P)}$.  In other words, each contour is an eigenvector of $P^{\otimes\sigma(P)}$.
Taking linear combinations of eigenvectors enables the construction of the eigenspace of $P^{\otimes\sigma(P)}$ (see Appendix~\ref{app:eigenspace}). This eigenspace is the set
of all possible periodic regimes (i.e., of all periods $\rho$,
$1\leq\rho\leq \sigma(P)$, obtained for all initial states
$\ve{x}(0)$), which corresponds to the set of all contour plots that
can be obtained.  

Theorem~\ref{thm:evector} (Appendix~\ref{app:eigenspace}) gives a method for constructing the eigenspace of an
irreducible matrix. Let $B=P^{\otimes\sigma(P)}$. Applying Theorem~\ref{thm:evector} to irreducible
$B$ will yield its eigenspace and consequently all possible contour
plots for the system $\ve{x}(k+1)=P\otimes\ve{x}(k)$. 

\begin{note}
If $P$ is irreducible, then it is not necessarily the case that $B=P^{\otimes\sigma(P)}$ will also be irreducible.  Nevertheless, in a highly connected lattice such as the regular network of this paper, it is more likely that $B$ is irreducible. If $B$ is reducible, then other methods to Theorem~\ref{thm:evector} must be employed.  One such method is Howard's policy improvement scheme (see \cite{Heid} or \cite{cochet1998}).
\end{note}

\begin{lemma}\label{lem:regimesnotunique}
Consider a max-plus system having irreducible timing dependency matrix $P$ with cyclicity $\sigma$.  For $\rho$ fixed
($1\leq\rho\leq\sigma$), period $\rho$ regimes are
not necessarily unique.
\begin{proof}
Each contour in the periodic regime is an eigenvector of $P^{\otimes \rho}$.  By taking linear combinations of eigenvectors, it is possible to construct linearly independent eigenvectors of $P^{\otimes \rho}$ such that the corresponding contours are also linearly independent.  Thus, there is no unique periodic regime of period $\rho$.
\end{proof}
\end{lemma}
The top of Figure~\ref{fig:syncandasynccontourplots}(a) shows the limiting contours of a size 20 system.  The contours depict a period 1 regime, therefore have the same shape.  Lemma~\ref{lem:regimesnotunique} says that, if $\sigma>1$ for that system, then other shapes of limiting contour may be possible for $\rho=1$ periodic regimes.  In other words, the limiting contour plot is not necessarily fixed, despite $\rho.\ve{\chi}$ being fixed.  The significance of this is that a
corresponding CA is asymptotically not unique, but dependent on the initial time $\ve{x}(0)$; this can affect the
time $x_i(k)$ of the $k^\text{th}$ update at node $i$ \emph{relative} to $x_j(k)$
(at node $j$), even when the `pattern' of \emph{consecutive} update times ($\ve{x}(k),\ve{x}(k+1),\ldots,\ve{x}(k+\rho),\ldots$) is
independent of $\ve{x}(0)$ (since the cycletime is independent of $\ve{x}(0)$).

%

\section{Cellular Automata in Max-Plus Time}\label{sec:OutlineCA}
We now present the first formalism for implementing a
cellular automaton asynchronously such that update times are determined by a max-plus
system.

\subsection{Contour Plot as a Foundation for Cellular Automata}\label{subsec:contourCA}
Let $s_i(k)$ denote the CA state of node $i$ at epoch $k$.  We are
concerned with Boolean CA states, so that $s_i(k)\in\{0,1\}$. The
unit $k$ is as used in the max-plus model which
updates the times $\ve{x}(k)$.  Thus, to be precise, $s_i(k)$ is the
CA state of node $i$ at time $x_i(k)$.  The \emph{CA state of
the system} is represented by the string $s_1(k)s_2(k)\cdots s_N(k)$, which can also be read
as the state of all nodes on contour $k$.  As a consequence, just as
we represented the vector $\ve{x}(k)$ by a contour, we can represent the
CA state $\ve{s}(k)$ by the same contour but with the addition that the coordinates $(i,x_i(k))$ now display the state $s_i(k)$ (e.g., in coloured form, where two different colours are used to distinguish the two states 1 and 0).

Recall the main events that are internal to node $i$; they occurring within
cycle $k$. These events are grouped in two: ``receive" and ``send".
Once the ``receive" CA states have all arrived, node $i$ applies a
CA rule on this set, to obtain the new state $s_i(k)$. If all
nodes have neighbourhood size $n$, the applied CA rule is the function $f:\{0,1\}^n\rightarrow\{0,1\}$ and the new state $s_i(k)$ is
calculated as
\begin{equation}
s_i(k) = f(s_{\mathcal{N}_i}(k-1)).
\end{equation}

\subsection{Cellular Automaton Space-Time Plot}
The classical one-dimensional CA is synchronous, so that the $k^\text{th}$
update time of each cell is the same.  Consequently, we can think of
such a system as having a contour plot that contains only horizontal
contours. Updates of the CA state of the system take place every one
time unit, thereby giving the synchronous CA a cycletime of 1.  The
time between contours in this system is thus of duration one,
although no such duration is depicted; for example, if $s_i(k)=1$ for
all $k$, then this is shown as a continuous vertical coloured block
in position $i$. 

\begin{figure}[!hbt]
\centering

 \includegraphics[scale=0.45,trim= 110 50 90 30,clip]{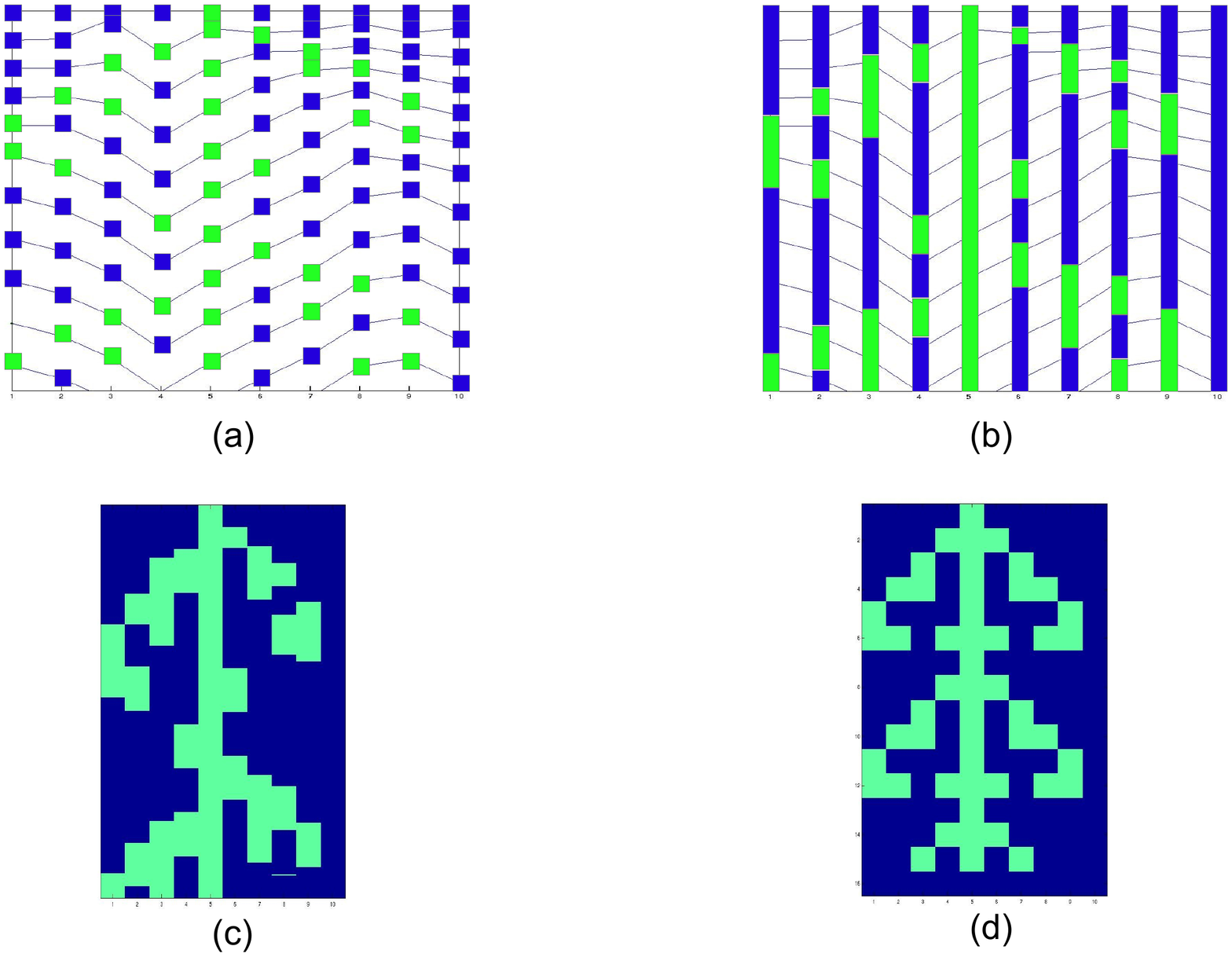}

\caption{Construction of the CA space-time plot in max-plus time.
The CA rule is ECA rule 150.  The initial CA state on contour 0 is
$s_5(0)=1$, $s_i(0)=0$ for $i\neq 5$.
 State 0 is coloured dark, state 1 is light.  In all figures, the vertical axis denotes real time, travelling down. (a) Contour plot with CA states indicated on each contour. (b) CA states indicated for all time by filling spaces between contours with memory. (c) Contours and space between nodes removed to obtain the CA space-time plot. (d) Classical (synchronous) CA space-time plot.}
 \label{fig:CAmaxplusevol}
\end{figure}

Despite varying contour shapes dictating the varying time gaps between
contours, it is simple and intuitive to construct the space-time plot for CA in max-plus time.  To illustrate, we use the example of a
regular $3$-nbhd network with size $N=10$. Let the positive
(diagonal) entries in matrix $A_{\xi}$ be represented by the vector $\ve{\xi}$ of processing times. We choose the entries in $\ve{\xi}$ at random with equal probability from all
integers between 1 and 30, whilst the non-zero entries in $T$ are
selected likewise from the integers between 1 and 10.  Taking the initial
time $\ve{x}(0)=\ve{u}$, we obtain a contour plot of update times by
iterating the max-plus system. We now address the CA state by assigning the depicted space between contours as \emph{memory}: for each
node, the CA state remains fixed until the time of update, which
corresponds to a contour.  For node $i$, the time $t_i$ that elapses between
contours implies that the storing of the CA state in memory can be represented as a vertical block of length $t_i$ (which is coloured
accordingly, depending on the CA state).  Correspondingly,  this may be depicted in a space-time plot,
the construction of which is shown in three stages in
Figure~\ref{fig:CAmaxplusevol}.

At this juncture, it is important to distinguish between the variables $t$ and $k$.  The term ``time" (or ``real time") now refers to a point $t\in\R^+$; it can be thought of as time as we know it.  $k$ maintains its role as a discrete epoch.  Thus, node $i$ carries a CA state for every point in real
time.  If we denote the state of node $i$ at real time $t$ as $s_i^{(t)}$, the CA state of node $i$ can now be understood in two ways: $s_i(k)$ denotes
the CA state on contour $k$, and it is discretely dynamic, whilst
$s_i^{(t)}$ represents the state in a dynamical system with a continuous
underlying real time $t$.  Thus, for example, on contour $k=2$, if $t=5.8$, we have $s_i(2)=s_i^{(5.8)}$.

Figure~\ref{fig:CAmaxplusevol}(d) is for comparison with
 Figure~\ref{fig:CAmaxplusevol}(c) and it shows the classical synchronous CA having the same initial time
$\ve{x}(0)$, initial CA state $\ve{s}(0)$ and CA rule. Whilst sharing
initial conditions and CA rule, the state $\ve{s}^{(t)}$ in both patterns will generally
differ. (This can be seen by simply drawing a horizontal line across
both patterns at time $t$ and reading off the state of each node at
that time).  The difference in pattern is clearly due to the
asynchrony of the max-plus system in
Figure~\ref{fig:CAmaxplusevol}(c). We show next
that we can characterise this difference somewhat and, in fact, map
the synchronous CA to the max-plus CA via the contour plot.

\subsection{Bijection}\label{subsec:bijection}
%
Recall the state transition graph.  In
Figure~\ref{fig:CAmaxplusevol}(d), the CA period is 6 and the
periodic orbit is the following set.
\begin{equation}
\begin{array}{rlllll}
\{ & 0000100000, & 0001110000, & 0010101000, & 0110101100, \\
 & 1000100010, & 1101110110 & \}. & & 
\end{array}
\end{equation}
If we consider the CA states on only the contours in
Figure~\ref{fig:CAmaxplusevol}(c) (which is seen better in
Figure~\ref{fig:CAmaxplusevol}(a)), we see that they are exactly the
same as the states in Figure~\ref{fig:CAmaxplusevol}(d). This is
a consequence of the max-plus model requiring \emph{all} neighbourhood
states to arrive before processing new CA states.  Indeed, this notion has also been mentioned (albeit briefly) in \cite[Page~1035]{WolNKS} under the heading of ``Intrinsic synchronization in cellular automata".

Given the same initial CA state $\ve{s}(0)$ and CA rule, let $\mathcal{S}$ and $\mathcal{M}$ denote the orbit of $\ve{s}(0)$ generated in the synchronous system and the max-plus system respectively.  Let $\ve{s}_{\mathcal{S}}(k)$ denote the CA state after $k$ iterations of the synchronous system; $\ve{s}_{\mathcal{M}}(k)$ denotes the CA state after $k$ iterations of the max-plus system.  The
model uses the same CA rule, applied to the same neighbourhoods, the
only difference being that the time of application of the rule is
different.  Then, after $k$ iterations of both systems, we clearly have $\ve{s}_{\mathcal{S}}(k)=\ve{s}_{\mathcal{M}}(k)$.  This defines a one-to-one and onto mapping - a bijection - between $\mathcal{S}$ and $\mathcal{M}$, and we say that both systems have the same state
transition graph (defined as the transitions between states on
contours).  Thus, the max-plus system need
not evolve the CA concurrently since the CA plot for the max-plus
system may be obtained from this mapping.   

In summary, the STG provides a deterministic form for predicting the
behaviour of the CA in max-plus time.  Each state $\ve{s}_{\mathcal{M}}(k)$ in the STG
does not necessarily correspond to a state in real time $t\in\R$ due to contours not necessarily being horizontal.  Nevertheless, the property of memory can be applied to ascertain such real time states.

\section{Conclusion}
We have shown that, when modelling discrete asynchronous systems, more attention needs to be paid to the internal processes of a cell.  This has resulted in the uncovering of a useful theory - that of max-plus algebra.  

Cellular automata are naturally well-suited to be modelled in max-plus time because this
model requires update on knowledge of all neighbours.  This has the additional benefit of cells updating only when they are ready.  Thus, whereas a fixed, periodic, global update time (as in the synchronous case) can be slower and less energy efficient, a max-plus asynchronous model consumes only the time and energy that local neighbourhoods require.  


The classical ECA corresponds to a strongly connected network in our model.  Thus, $P$ is irreducible, and this ensures periodic
behaviour, which is not what is usually associated with the word ``asynchrony".  
We comment that periodic behaviour is expected also when $P$ is reducible; in this case, the theory is similar to what we have covered here but the resultant cycletime vector is not necessarily uniform \cite{Heid}.  We have therefore seen that not only is it an efficient system for the timing of asynchronous CA, the max-plus system is also a simple, deterministic asynchronous model.

We can get a visual sense of such periodicity via the contour plot, which couples the asynchronous update times with the CA.  The eigenspace of $P^{\otimes{\sigma(P)}}$ yields the range of contour plots that can be generated due to $P$; this enables us to get a sense of the range of the corresponding CA space-time patterns that are possible.  By identifying a bijection between synchronous and max-plus CA, we can further narrow the aforementioned range of CA patterns since CA states on contours can be obtained directly from the synchronous CA.  We note that this notion of memory has been suggested previously (see \cite{fates2013}, \cite{nehaniv},  \cite{WolNKS} and the references therein); this bijection is captured under the topic of ``causal invariance" in \cite{WolNKS}, wherein the the same causal network (i.e., structure of events - past and present) emerges, irrespective of CA states.  In this paper, however, we have identified a numerical link between the idea and the theory in the form of max-plus algebra.

Importantly, the contour plot also allows one to see that the system visits many more `interim' states (in memory) in real time - these are transient CA states that illustrate the local dynamics.  Thus, although the STG is the same as the synchronous STG, max-plus CA gives more information: we can now visualise exactly \emph{how} one state in the STG evolves into another.  Max-plus algebra offers promising scope for assigning numerical measures to the states in memory; the parameters $\xi_i$, $\tau_{ij}$, along with the cycletime, are likely to play major roles for this purpose.  Such classification scales up to the CA space-time pattern as a whole, and attempts to do this are the focus of further work.  

Synchronous CA have previously been held up as models for patterns seen in nature (such as those seen on seashells and the growth of snowflakes).  Certain probabilistic CA improved on this to account for random fluctuations in the growth processes.  We hypothesise that a max-plus algebraic approach adds further realism because it considers the actual processing and delay times that may be present within the chemical reactions and relies not on probability.  One might conjecture the aforementioned transient states in memory to be precisely the processes observed during the natural construction of such patterns as found on seashells and snowflakes, that is, apparent deviations from `normal' growth may, in fact, be part of a transient phase as opposed to some random fault. Noting the absence of significant literature on the matter, we believe Boolean networks provide a notable avenue to exploit; it would be interesting to see how the various topologies of these networks behave under max-plus time.  How does a max-plus model impact on applications such as Kauffman's genetic regulatory network in \cite{Kauf1}?  Moreover, in light of connections to the idea of causal invariance in \cite{WolNKS}, might max-plus algebra now provide a numerical measure for the related topics (e.g., of space, time and relativity) contained therein?

In addition to the above work, this paper has also laid the theoretical groundwork for extending the asynchronous model itself.  For instance, one such extension includes the minimum operator, and can subsequently be used to conduct further studies of asynchronous CA, particularly to better describe an intended application.  Preliminary results of these new ``max-min-plus" models are less predictable, particularly manifested by an absence of the bijection described in Section~\ref{subsec:bijection} \cite{PatelMaxmin}.

\section*{Acknowledgments}
The novel idea in this paper was conceived by this coauthor but ill health - and his subsequent passing in July 2014 - prevented involvement in its preparation.

This work was initiated as part of the Centre for Interdisciplinary Computational and Dynamical Analysis (CICADA) project which was funded by the Engineering and Physical Sciences Research Council (EPSRC).  E. L. Patel is now funded by the EPSRC Horizon: Digital Economy project (EP/G065802/1).


\appendix

\section{Maximal strongly connected subgraph}\label{app:MSCS}
We say that node $i$ \emph{communicates with} node $j$, denoted $i\leftrightarrow j$, if either $i=j$ or $i\rightarrow j$ and $j\rightarrow i$.  Note that we allow a solitary node to communicate with itself, even if there is no self-loop attached to it.  

It is, thus, possible to partition the node set $V$ of a graph into disjoint subsets
$V_i$ such that $V=V_1\cup V_2\cup\cdots\cup V_q$, where each subset
$V_i$ contains nodes that communicate with each
other but not with other nodes of $V$.  By taking $V_i$ together
with arc set $E_i$, each of whose arcs has start node and end node
in $V_i$, we obtain the subgraph $\mathcal{G}_i=(V_i,E_i)$.  We
call this subgraph a \emph{maximal strongly connected subgraph}
(MSCS) of $\mathcal{G}=(V,E)$.

\section{Asymptotic behaviour of the max-plus system} \label{app:asymptotic}

\begin{theorem}\label{thm:evalueirr}
Let $A\in\R_{\max}^{n\times n}$ be irreducible.  Then $A$ possesses
a unique eigenvalue, denoted $\lambda(A)$, which is finite
($\neq\varepsilon$).  Moreover, this eigenvalue is equal to the
maximal average weight of elementary circuits in $\mathcal{G}(A)$.
Let $c$ denote an elementary circuit of $\mathcal{G}(A)$.  Denote
the set of all elementary circuits of $\mathcal{G}(A)$ by
$\mathcal{C}(A)$. Then
\begin{equation}
\lambda(A)=\max_{c\in\mathcal{C}(A)}\frac{|c|_w}{|c|_l}.
\end{equation}
\begin{proof}
See \cite[Theorem 2.9]{Heid}.
\end{proof}
\end{theorem}

The following theorem shows that the cycletime vector $$\ve{\chi}=\left(\lim_{k\rightarrow\infty}\frac{x_1(k)}{k},\lim_{k\rightarrow\infty}\frac{x_2(k)}{k},\ldots,\lim_{k\rightarrow\infty}\frac{x_N(k)}{k}\right)^\top$$ is unique.
\begin{theorem} \label{thm:cycletimeindependent}
Consider the recurrence relation $\ve{x}(k+1)=A\otimes\ve{x}(k)$ for
$k\geq0$ and $A\in\R_{\max}^{n\times n}$ irreducible.  For some
$\ve{x}_\star(0)\in\R_{\max}^n$ whose elements are all finite, if the limit $\lim_{k\rightarrow\infty}\frac{A^k\otimes\ve{x}_\star(0)}{k}$ exists, then
this limit is the same for any initial condition
$\ve{x}(0)\in\R_{\max}^n$ whose elements are all finite.
\begin{proof}
See \cite[Theorem 3.11]{Heid}.
\end{proof}
\end{theorem}
In fact, the condition of irreducibility can be relaxed to that of reducibility, as long as each node has at least one predecessor node in  $\mathcal{G}(A)$; this corresponds to all rows of $A$ containing
at least one non-zero ($\neq-\infty$) element.  Thus, in the latter (relaxed) case, we obtain a
cycletime vector whose elements may not necessarily be identical. However, for $A$
irreducible, each element of $\ve{\chi}$ turns out to be the same - specifically the eigenvalue of $A$ - as stated in the following lemma.
\begin{lemma}\label{lem:allnodessamecycle}
For the recurrence relation $\ve{x}(k+1)=A\otimes\ve{x}(k)$ with
$k\geq0$, let $A\in\R_{\max}^{n\times n}$ be an irreducible matrix
having eigenvalue $\lambda\in\R$.  Then, for $i=1,2,\ldots,n$,
\begin{displaymath}
\lim_{k\rightarrow\infty}\frac{x_i(k)}{k}=\lambda
\end{displaymath}
for any initial condition $\ve{x}(0)\in\R^n$.
\begin{proof}
See \cite[Lemma 3.12]{Heid}.
\end{proof}
\end{lemma}

\section{Eigenspace of an irreducible matrix}\label{app:eigenspace}
As in conventional linear algebra, eigenvectors
are not unique in max-plus algebra because they are
defined up to scalar multiplication.  (It can easily be shown that, if $\ve{v}$ and $\ve{w}$ are eigenvectors of $A\in\R_{\max}^{n\times n}$ associated
with eigenvalue $\lambda$, then, for $\alpha,\beta\in\R_{\max}$,
$\alpha\otimes\ve{v}\oplus\beta\otimes\ve{w}$ is also an
eigenvector of $A$.)

Consider the definition of the \emph{Kleene
star} for any $A\in\R_{\max}^{n\times n}$:
\begin{equation}
A^*\stackrel{\text{def}}{=}\bigoplus_{k=0}^{\infty}A^{\otimes k}.
\end{equation}
It is known that, if circuit weights
in $\mathcal{G}(A)$ are nonpositive, then the Kleene star of a
square matrix over $\R_{\max}$ exists \cite{Heid}.  Denote the critical graph of $A$ as
$\mathcal{G}^{cr}(A)=(V^{cr}(A),E^{cr}(A))$ and the normalised
matrix $\hat{A}=-\lambda\otimes A$.  
\begin{theorem}\label{thm:evector}
Let $A\in\R_{\max}^{n\times n}$ be irreducible and consider
$\hat{A}^*$ to be the Kleene star of $\hat{A}=-\lambda\otimes A$.
\begin{enumerate}
\item If node $i$ belongs to $\mathcal{G}^{cr}(A)$, then
$[\hat{A}^*]_{\cdot i}$ is an eigenvector of $A$.
\item The eigenspace of $A$ is
\begin{equation}
V(A)=\left\{\ve{v}\in\R_{\max}^n|\ve{v}=\bigoplus_{i\in
V^{cr}(A)}a_i\otimes[\hat{A}^*]_{\cdot i} \quad\text{for
}a_i\in\R_{\max}\right\}.
\end{equation}
\item For $i,j$ belonging to $\mathcal{G}^{cr}(A)$, there exists
$a\in\R$ such that
\begin{equation}
a\otimes[\hat{A}^*]_{\cdot i}=[\hat{A}^*]_{\cdot j}
\end{equation}
if and only if $i$ and $j$ belong to the same MSCS of
$\mathcal{G}^{cr}(A)$.
\end{enumerate}
\begin{proof}
In \cite[Thoerem 4.5]{Heid}.
\end{proof}
\end{theorem}

\section{Hasse diagram}\label{app:hasse}
With elements taken from \cite{Cam}, we define the Hasse diagram formally.  We first require the ``happened before" relation, denoted by $\prec$.

\begin{defn}\label{def:happbefore}
The relation “$\prec$” on a set of events is defined by the following
conditions:
\begin{enumerate}
\item If the events $a$ and $b$ are processed by the same processor, and $a$ occurs before $b$, then $a\prec b$.
\item If $a$ is the sending of a message by processor $A$ and $b$ is the receipt of the
message by another processor $B$, then $a\prec b$.
\item If $a \prec b$, and $b\prec c$, then
$a \prec c$.
\end{enumerate}
\end{defn}
We say that two distinct events $a$ and $b$ are “concurrent” if $a\nprec b$ and $b \nprec a$.
We also assume the properties of irreflexivity, that is, $a\nprec  a$, and antisymmetry on
the times of events, that is, if $a_t$ and $b_t$ represent the times of events $a$ and $b$, then
$a_t \nprec b_t$, $b_t \nprec a_t \Rightarrow a_t = b_t$.

In fact, Definition~\ref{def:happbefore}, along with the properties of irreflexivity and antisymmetry, defines a \emph{partial ordering} on the set $X$ of all events in our system.  We say that the set $X$, along with the relation $\prec$,
forms a partially ordered set (``poset").  In the following, the pair $(X,\prec)$ denotes a poset. 
\begin{defn}
Let $x$ and $y$ be distinct elements of a poset $(X,\prec)$.  $y$ is
said to \emph{cover} $x$ if $x\prec y$ but no element $z$ satisfies
$x\prec z\prec y$.
\end{defn}
\begin{defn}\label{def:antichain}
Let the set $X_1$ of $n$ elements $\{x_1,x_2,\ldots,x_n\}$ be a
subset of $(X,\prec)$ such that each element may be
totally ordered according to $\prec$ as $x_1\prec x_2\prec\cdots
\prec x_n$. Then $X_1$ is a \emph{chain}.  The subset $X_2\in X$ is
called an \emph{antichain} if and only if no elements of $X_2$ may
be totally or partially ordered.
\end{defn}
These definitions now allow us to formulate the following definition.
\begin{defn}[Hasse diagram] \label{def:Hasse}
The \emph{Hasse diagram} of a poset $(X,\prec)$ is a graph drawn in
the Euclidean plane such that each element of the poset is
represented by a unique vertex in the graph. Each covering pair
$x\prec y$, is depicted by a directed arc from $x$ to $y$, where the point representing $x$ is below
the point representing $y$ (i.e., it has smaller $Y$-coordinate).
\end{defn}

\end{document}